\documentclass[11pt]{amsart}
\hoffset=- 2cm \voffset=- 1cm 
\numberwithin{equation}{section}
\usepackage[latin1]{inputenc}
\usepackage{amsmath,amsfonts,eufrak,amssymb,graphics}
\usepackage[mathcal]{eucal}

\def\eq#1{(\ref{#1})}
\def\neweq#1{\begin{equation}\label{#1}}
\def\endeq{\end{equation}}

\def\ep{\varepsilon}

\def\ph{\varphi}
\def\RR{{\bf R}}

\def\di{\displaystyle}
\def\ri{\rightarrow}

\def\ov{\overline}

\def\neweq{\begin{equation}}
\def\endeq{\end{equation}}

\def\di{\displaystyle}
\def\ri{\rightarrow}
\def\ov{\overline}

\def\ep{\varepsilon}

\def\RR{{\mathbf R}}

\theoremstyle{definition}
\newtheorem{definition}{Definition}[section]
\theoremstyle{remark}
\newtheorem{remark}{Remark}[section]
\newtheorem{example}{Example}[section]

\theoremstyle{plain}
\newtheorem{theorem}{Theorem}[section]
\newtheorem{lemma}{Lemma}[section]
\newtheorem{proposition}{Proposition}[section]

\begin{document}
\title[Nonlinear problems with boundary blow-up]
{Nonlinear problems with boundary blow-up: a Karamata regular
variation theory approach}
\author[F. C\^{i}rstea]{Florica Corina C\^{\i}rstea}
\address{F. C\^{i}rstea: Centre for Mathematics and its Applications,
Mathematical Sciences Institute, The Australian National University,
Canberra, ACT 0200, 
Australia}
\email{Florica.Cirstea@maths.anu.edu.au}

\author[V. R\u adulescu]{Vicen\c tiu R\u adulescu}
\address{V. R\u adulescu: Department
of Mathematics, University of Craiova, 13 A. I. Cuza Street,
200585 Craiova, Romania.} \email{vicentiu.radulescu@ucv.ro}
\urladdr{http://inf.ucv.ro/\textasciitilde radulescu}
\date{March 2005}
\thanks{F. C\^ {\i}rstea was supported by the Australian Government
through DETYA and Victoria University (Melbourne) under the IPRS Programme}
\begin{abstract}
We study the uniqueness and expansion properties of the positive
solution of the logistic equation $\Delta u+au=b(x)f(u)$ in a
smooth bounded domain $\Omega$, subject to the singular boundary
condition $u=+\infty$ on $\partial\Omega$. The absorption term $f$
is a positive function satisfying the Keller--Osserman condition
and such that the mapping $f(u)/u$ is increasing on $(0,+\infty)$.
We assume that $b$ is non-negative, while the values of the real
parameter $a$ are related to an appropriate semilinear eigenvalue
problem. Our analysis is based on the Karamata regular variation
theory.
\end{abstract}
\maketitle

\section{Introduction and main results}

Let $\Omega\subset \RR^N$ $(N\geq 3)$ be a smooth bounded domain.

Consider the semilinear elliptic equation
\neweq\Delta u+au=b(x)f(u)\qquad \mbox{in}\ \Omega, \label{0a}
\endeq
where $f\in C^1[0,\infty)$, $a\in \RR$ is a parameter and $b\in
C^{0,\mu}(\overline{\Omega})$ satisfies $b\geq 0$, $b\not\equiv 0$
in $\Omega$. Such equations are also known as the stationary
version of the Fisher equation \cite{fi} and the
Kolmogoroff--Petrovsky--Piscounoff  equation \cite{KPP}  and they
have been studied by Kazdan--Warner \cite{kw}, Ouyang \cite{ou},
del Pino \cite{dP} and Du--Huang \cite{dh}.

Note that if $f(u)=u^{(N+2)/(N-2)}$, then (\ref{0a}) originates
from the Yamabe problem, which is a basic problem in Riemannian
geometry (see, e.g., \cite{lee}).

The existence of positive solutions of (\ref{0a}) subject to the
Dirichlet boundary condition, $u=0$ on $\partial\Omega$, has been
intensively studied in the case $f(u)=u^p$, $p>1$ (see
\cite{al_t1}, \cite{ag}, \cite{Da}, \cite{dP}, \cite{fklm} and
\cite{ou}); this problem is a basic population model (see
\cite{hs}) and it is also related to some prescribed curvature
problems in Riemannian geometry (see \cite{kw} and \cite{ou}).
Moreover, if $b>0$ in $\ov{\Omega}$, then it is referred to as the
logistic equation and it has a unique positive solution if and
only if $a>\lambda_1(\Omega)$, where $\lambda_1(\Omega)$ denotes
the first eigenvalue of $(-\Delta)$ in $H_0^1(\Omega)$.

In the understanding of (\ref{0a}) an important role is played by
the interior of the zero set of $b$: \[ \Omega_0:={\rm int}\,
\{x\in \Omega:\ b(x)=0\}.\]

We assume, throughout this paper, that $\Omega_0$ is connected
(possibly empty), $\ov{\Omega}_0\subset \Omega$ and $b>0$ in
$\Omega\setminus \ov{\Omega}_0$. Note that we allow $b\geq 0$ on
$\partial\Omega$. Let $\partial\Omega_0$ satisfy an exterior cone
condition and $\lambda_{\infty,1}$ be the first Dirichlet
eigenvalue of $(-\Delta)$ in $H_0^1(\Omega_0)$ (with
$\lambda_{\infty,1}=\infty$ if $\Omega_0=\emptyset$).

By a \emph{large} (or \emph{blow-up}) solution of (\ref{0a}), we
mean any non-negative $C^2(\Omega)$-solution of (\ref{0a}) such
that $u(x)\to \infty$ as $d(x):={\rm dist}\,(x,\partial\Omega)\to
0$.

Assuming that $f$ satisfies
\begin{equation} \tag{$A_1$}
f\in C^1[0,\infty)\ \mbox{is non-negative and } f(u)/u\ \mbox{is
increasing on } (0,\infty),
\end{equation}
then, necessarily $f(0)=0$, and by the strong maximum principle,
any non-negative classical solution of (\ref{0a}) is positive in
$\Omega$ unless it is identically zero. Consequently, any large
solution of (\ref{0a}) is positive. Moreover, it is well known
(see, e.g., Remark~1.1 in \cite{cr2}) that in this situation, the
Keller--Osserman condition
\begin{equation} \tag{$A_2$}
\int_{1}^{\infty}\frac{dt}{\sqrt{F(t)}}< \infty\,, \ \
\mbox{where}\ \ F(t)=\di\int_{0}^{t} f(s)\,ds \end{equation} is
necessary for the existence of large solutions of (\ref{0a}).

When $(A_1)$ and $(A_2)$ hold, Theorem 1.1 in \cite{cr2} shows
that (\ref{0a}) possesses large solutions if and only if $a<
\lambda_{\infty,1}$. The hypothesis $(A_1)$ is inspired by
\cite{al_t1}, where it is developed an exhaustive study of
positive solutions of (\ref{0a}), subject to $u=0$ on
$\partial\Omega$.

Our major goal is to advance innovative methods to study the
uniqueness and asymptotic behavior of large solutions of
(\ref{0a}). We develop the research line opened up in \cite{cr3}
to gain insight into the two-term asymptotic expansion of the
large solution near $\partial\Omega$. Our approach relies
essentially on the \emph{regular variation theory} (see \cite{bgt}
and section~\ref{sec2}) not only in the statement but in the proof
as well. This enables us to obtain significant information about
the qualitative behavior of the large solution to (\ref{0a}) in a
general framework that removes previous restrictions in the
literature.

We point out that, despite a long history and intense research on
the large solutions, the regular variation theory arising in
probability theory has not been exploited before in this context.

Singular value problems having large solutions have been initially
studied for the special case $f(u)=e^u$ by Bieberbach \cite{bie}
(if $N=2$). Problems of this type arise in Riemannian geometry.
More precisely, if a Riemannian metric of the form
$|ds|^2=e^{2u(x)}|dx|^2$  has constant Gaussian curvature $-g^2$
then $\Delta u=g^2e^{2u}$. This study was continued by Rademacher
\cite{rad} (if $N=3$) in connection with some concrete questions
arising in the theory of Riemann surfaces, automorphic functions
and in the theory of the electric potential in a glowing hollow
metal body.

The question of large solutions was later considered in
$N$-dimensional domains and for other classes of nonlinearities
(see \cite{be}, \cite{bm}--\cite{bm2}, \cite{cr0}--\cite{cr4}, \cite{de2},
\cite{dh}, \cite{gls}, \cite{1}, \cite{kon}--\cite{3},  \cite{8}, \cite{mv1}--\cite{mv2}, 
\cite{matero}, \cite{2}, \cite{rrv}).

In higher dimensions the notion of Gaussian curvature has to be
replaced by the scalar curvature. It turns out that if a metric of
the form $|ds|^2=u(x)^{4/(N-2)}|dx|^2$ has constant scalar
curvature $-g^2$, then $u$ satisfies (\ref{0a}) for
$f(u)=u^{(N+2)/(N-2)}$, $a=0$ and $b(x)=[(N-2)g^2]/[4(N-1)]$. In a
celebrated paper, Loewner and Nirenberg \cite{8} described the
precise asymptotic behavior at the boundary of large solutions to
this equation  and used this result in order to establish the
uniqueness of the solution. Their main result is derived under the
assumption that $\partial\Omega$ consists of the disjoint union of
finitely compact $C^\infty$ manifolds, each having codimension
less than $N/2+1$. More precisely, the uniqueness of a large
solution is a consequence of the fact that every large solution
$u$ satisfies
\neweq u(x)={\mathcal   E}(d(x))+o({\mathcal   E}(d(x)))
\quad \mbox{as}\ d(x)\to 0,\label{0ln}
\endeq
where ${\mathcal   E}$ is defined by
\neweq \int_{{\mathcal   E}(t)}^\infty \frac{ds}{\sqrt{2F(s)}}=
\left(\frac{(N-2)g^2}{4(N-1)}\right)^{1/2}\,t,\quad \mbox{for all}
\ t>0. \label{as}
\endeq

Kondrat'ev and Nikishkin \cite{kon} established the uniqueness of
a large solution for the case $a=0$, $b=1$ and $f(u)=u^p$ ($p\geq
3$), when $\partial\Omega$ is a $C^2$-manifold and $\Delta$ is
replaced by a more general second order elliptic operator.

Dynkin \cite{dyn} showed that there exist certain relations
between hitting probabilities for some Markov processes called
superdiffusions and maximal solutions of (\ref{0a}) with $a=0$,
$b=1$ and $f(u)=u^p$ ($1<p\leq 2$). By means of a probabilistic
representation, a uniqueness result in domains with non-smooth
boundary was established by le Gall \cite{gall} when $p=2$. We
point out that the case $p=2$ arises in the study of the subsonic
motion of a gas. In this connection the question of uniqueness is
of special interest.

Recently, \cite{gls} gives the uniqueness and exact two-term
asymptotic expansion of the large solution of (\ref{0a}) in the
special case $f(u)=u^p$ ($p>1$), $b>0$ in $\Omega$ and $b\equiv 0$
on $\partial\Omega$ such that
\neweq
b(x)=C_0 [d(x)]^{\gamma}+o([d(x)]^{\gamma})\ \mbox{as}\ d(x)\to
0,\ \mbox{for some constants}\ C_0,\gamma>0. \label{ass}
\endeq
It was shown there that the degenerate case $b\equiv 0$ on
$\partial \Omega$ is a \emph{natural} restriction for $b$
inherited from the logistic equation.

\medskip To present our main results, we briefly recall some notions
from Karamata's theory (see \cite{bgt} or \cite{se}); more details
are provided in section~\ref{sec2}.

A positive measurable function $R$ defined on $[A, \infty)$, for
some $A>0$, is called {\em regularly varying  with index} $q\in
\RR$, written $R\in RV_q$, provided that \[ \lim_{u\to
\infty}\frac{R(\lambda u)}{R(u)}=\lambda^q,\qquad \mbox{for all}\
\lambda>0.\] When the index $q$ is zero, we say that the function
is {\em slowly varying}.

Clearly, if $R\in RV_q$, then $L(u):=R(u)/u^q$ is a slowly varying
function.

Let ${\mathcal  K}$ denote the set of all positive, non-decreasing
$k\in C^1(0,\nu)$ that satisfy
\[\lim_{t\searrow 0}
\left(\frac{\int_0^t k(s)\,ds}{k(t)}\right):=\ell_0\ \ \mbox{and}
\ \ \lim_{t\searrow 0}\left(\frac{\int_0^t
k(s)\,ds}{k(t)}\right)':=\ell_1. \] Notice that $\ell_0=0$ and
$\ell_1\in [0,1]$, for every $k\in {\mathcal K}$. Thus, ${\mathcal
K}={\mathcal K}_{(01]}\cup {\mathcal K}_0$, where
\[{\mathcal K}_{(01]}=\{k\in {\mathcal K}:\ \ 0<\ell_1\leq 1\}\quad
\mbox{and} \quad {\mathcal K}_0=\{k\in {\mathcal K}:\ \
\ell_1=0\}.\] The exact characterization of ${\mathcal K}_{(01]}$
and ${\mathcal K}_0$ will be provided in section~\ref{gyr}.

If $H$ is a non-decreasing function on $\RR$, then we define the
(left continuous) inverse of $H$ by \[ H^\leftarrow (y)=\inf\{s:\
H(s)\geq y\}.\]

Our first result establishes the uniqueness of the large solution
of (\ref{0a}).

\begin{theorem}\label{uni1}
Let $(A_1)$ hold and $f\in RV_{\rho+1}$ with $\rho>0$. Suppose
there exists $k\in {\mathcal   K}$ such that
\begin{equation} \label{B} b(x)=k^2(d)+o(k^2(d))\ \mbox{as } d(x)\to
0.
\end{equation}
Then, for any $a\in (-\infty,\lambda_{\infty,1})$, $(\ref{0a})$
admits a unique large solution $u_a$. Moreover, the asymptotic
behavior is given by
\neweq u_a(x)=[2(2+\ell_1 \rho)/\rho^2]^{1/\rho}\, \ph(d)+o(\ph(d)) \quad \mbox{as } d(x)\to 0,
\label{comp} \endeq where $\ph$ is defined by
\neweq \frac{f(\ph(t))}{\ph(t)}=\frac{1}{\left(\int_0^t
k(s)\,ds \right)^2},\quad \mbox{for } t>0\ \mbox{small}.
\label{defp}
\endeq
\end{theorem}

Under the assumptions of Theorem~\ref{uni1}, let $r(t)$ satisfy
$\lim_{t\searrow 0}\left(\int_0^t k(s)\,ds\right)^2 r(t)=1$ and
$\widehat{f}(u)$ be chosen such that $\lim_{u\to
\infty}\widehat{f}(u)/f(u)=1$ and $j(u)=\widehat{f}(u)/u$ is
non-decreasing for $u>0$ large. Then, $\lim_{t\searrow 0} \ph(t)/
\widehat{\ph}(t)=1$, where $\ph$ is defined by (\ref{defp}) and $
\widehat{\ph}(t)=j^\leftarrow (r(t))$ for $t>0$ small.

The behavior of $\ph(t)$ for small $t>0$ will be described in
section~\ref{gyr}. In particular, if $k\in {\mathcal K}$ with
$\ell_1\not=0$, then $\ph(1/u)\in RV_{2/(\rho\ell_1)}$. In
contrast, if $k\in {\mathcal K}$ with $\ell_1=0$, then
$\ph(1/u)\not\in RV_q$, for all $q\in \RR$ (see Remark~\ref{see}).

\begin{remark} Theorem~\ref{uni1} improves the main result in \cite{cr3}, where
assuming that $f' \in RV_{\rho}$ (which yields $f\in
RV_{\rho+1}$), we prove \begin{equation} \label{fjk} u_a(x)=\xi_0
h(d)+o(h(d)) \quad \mbox{as } d(x)\to 0,\end{equation} where
$\xi_0=\left(\frac{2+\ell_1\rho}{2+\rho}\right)^{1/\rho}$ and $h$
is given by
\neweq \int_{h(t)}^\infty \frac{ds}{\sqrt{2F(s)}}=\int_0^t
k(s)\,ds,\quad \mbox{for } t>0 \mbox{ small}. \label{defi} \endeq
\end{remark}

\begin{remark} Theorem~\ref{uni1} recovers the uniqueness results of \cite{8} and
\cite{gls}. Note that for $k(t)=[(N-2)g^2/4(N-1)]^{1/2}$ in \eq{B}
and $f(u)=u^{(N+2)/(N-2)}$, (\ref{comp}) reduces to relation
(\ref{0ln}), prescribed by Loewner and Nirenberg \cite{8} for
their problem. Moreover, if $f(u)=u^p$ (with $p=\rho+1>1$) and
$k(t)=\sqrt{C_0}\,t^{\gamma/2}$ ($C_0,\gamma>0$), then we regain
the uniqueness result of \cite{gls}.
\end{remark}

The next objective is to find the two-term blow-up rate of $u_a$
when \eq{B} is replaced by
\begin{equation} \label{tilde B} b(x)=k^2(d)(1+\widetilde c
d^\theta+o(d^\theta))\ \mbox{as } d(x)\to 0,\end{equation} where
$\theta>0$, $\widetilde c\in \RR$ are constants. To simplify the
exposition, we assume that $f' \in RV_\rho$ $(\rho>0)$, which is
equivalent to $f(u)$ being of the form
\neweq
f(u)=Cu^{\rho+1}{\rm exp}\left\{\int_B^u
\frac{\phi(t)}{t}\,dt\right\}, \quad \forall u\geq B, \label{form}
\endeq for some constants $B,\,C>0$, where $\phi\in C[B,\infty)$
satisfies $\lim_{u\to \infty} \phi(u)=0$. In this case, $f(u)/u$
is increasing on $[B,\infty)$ provided that $B$ is large enough.

We prove that the two-term asymptotic expansion of $u_a$ near
$\partial\Omega$ depends on the chosen subclass for $k\in
{\mathcal  K}$ and the additional hypotheses on $f$ (by means of
$\phi$ in (\ref{form})).

Let $-\rho-2<\eta\leq 0$ and $\tau, \zeta>0$. We define
\[\begin{aligned}
{\mathcal  F}_{\rho\eta}&=\left\{f'\in RV_{\rho}\ (\rho>0):\
\mbox{either}\ \phi\in RV_\eta\ \mbox{or}\
-\phi\in RV_\eta\right\},\\
{\mathcal  F}_{\rho0,\tau}&= \{f\in {\mathcal  F}_{\rho 0}:\
\lim_{u\to \infty}(\ln u)^\tau \phi(u)=\ell^\star\in
\RR\},\\
{\mathcal  K}_{(01],\tau} &= \left\{k\in {\mathcal  K}_{(01]}:\
\lim_{t\searrow 0}(-\ln t)^\tau\left[\left(\frac{\int_0^t
k(s)\,ds}{k(t)}\right)'-\ell_1\right]:=L_\sharp\in \RR\right\},\\
{\mathcal  K}_{0,\zeta} &=\left\{ k\in {\mathcal  K}_0:\
\lim_{t\searrow 0}\frac{1}{t^\zeta}\,\left(\frac {\int_0^t
k(s)\,ds}{k(t)}\right)':=L_\star\in \RR\right\}.
\end{aligned} \]
Further in the paper, $\eta$, $\tau$ and $\zeta$ are understood in
the above range.

For the sake of comparison, we state here the following result.

\begin{theorem}\label{uni2}
Suppose $(A_1)$, {\rm \eq{tilde B}} with $k\in {\mathcal
K}_{0,\zeta}$, and one of the following growth conditions at
infinity:
\begin{enumerate}
\item[{\rm (i)}] $f(u)=Cu^{\rho+1}$ in a neighborhood of infinity
(i.e., $\phi\equiv 0$ in $(\ref{form})$); \item[{\rm (ii)}] $f\in
{\mathcal  F}_{\rho\eta}$ with $\eta\not=0$; \item[{\rm (iii)}]
$f\in {\mathcal  F}_{\rho0,\tau_1}$ with $\tau_1=\varpi/\zeta$,
where $\varpi=\min\{\theta,\zeta\}$. \end{enumerate} Then, for any
$a\in (-\infty,\lambda_{\infty,1})$, the two-term blow-up rate of
$ u_a$ is
\begin{equation}
u_a(x) =\xi_0 h(d)(1+\chi d^\varpi+o(d^\varpi))\quad \mbox{as }
d(x)\searrow 0 \label{chi1} \end{equation} where $h$ is given by
{\rm \eq{defi}}, $\xi_0=[2/(2+\rho)]^{1/\rho}$ and
\[
\chi=\left\{\begin{aligned} & \frac{L_\star}{2} {\rm
Heaviside(\theta-\zeta)}-\frac{\widetilde c}{\rho}{\rm
Heaviside(\zeta-\theta)}:=\chi_1\ \mbox{if}\ \ {\rm (i)} \
\mbox{or}\
{\rm (ii)} \ \mbox{holds},\\
& \chi_1- \frac{\ell^\star}{\rho} \left[ \frac{\rho \zeta
L_\star}{2(1+\zeta)}\right]^{\tau_1} \left( \frac{1}{\rho+2}+\ln
\xi_0 \right)\ \ \mbox{if }f \mbox{ obeys } {\rm (iii)}.
\end{aligned} \right.
\]
\end{theorem}

Theorem~\ref{uni2} is a consequence of \cite[Theorem~1]{cr4} and
Proposition~\ref{lok}.

\begin{theorem}\label{uni3}
Suppose $(A_1)$, {\rm \eq{tilde B}} with $k\in {\mathcal
K}_{(01],\tau}$, and one of the following conditions:
\begin{enumerate}
\item[{\rm (i)}] $f\in {\mathcal  F}_{\rho\eta}$ with $\eta
L_\sharp\not=0$; \item[{\rm (ii)}] $f\in {\mathcal
F}_{\rho0,\tau}$ with $[\ell^\star(\ell_1-1)]^2+L_\sharp^2\not=0$.
\end{enumerate}
Then, for any $a\in (-\infty,\lambda_{\infty,1})$, the two-term
blow-up rate of $ u_a$ is
\begin{equation}
u_a(x)=\xi_0 h(d)[1+\widetilde \chi\,(-\ln d)^{-\tau}+ o((-\ln
d)^{-\tau})]\quad \mbox{as } d(x)\searrow
0,\label{chi3}\end{equation} where $h$ is given by {\rm
\eq{defi}}, $\xi_0=[(2+\ell_1\rho)/(2+\rho)]^{1/\rho}$ and \neweq
\label{hii} \widetilde \chi=\left\{\begin{aligned} &
\frac{L_\sharp}{2+\rho\ell_1}:=\chi_2
\ \ \mbox{if}\ {\rm (i)}\ \mbox{holds},\\
& \chi_2- \frac{\ell^\star}{\rho} \left( \frac{ \rho \ell_1}{2}
\right)^\tau \left[ \frac{2(1-\ell_1)}{(\rho+2)
(\rho\ell_1+2)}+\ln \xi_0\right]\ \ \mbox{if } f \mbox{ obeys}\
{\rm (ii)}.
\end{aligned} \right.\endeq
\end{theorem}

\begin{remark}
Note that Theorems~\ref{uni2} and \ref{uni3} distinguish from
Theorem~1 in \cite{gls}, which treats the particular case
$f(u)=u^p$ $(p>1)$, $\Omega_0=\emptyset$, $k(t)=\sqrt{C_0
t^{\gamma}}$ ($C_0,\gamma>0$) and $\theta=1$ in (\ref{tilde B}).
The second term in the asymptotic expansion of $u_a$ near
$\partial\Omega$ involves in \cite{gls} both the distance function
$d(x)$ and the mean curvature of $\partial\Omega$.

Theorem~\ref{uni2} admits the case $f(u)=u^p$ assuming that $k\in
{\mathcal K}_{0,\zeta}$, while the alternative (ii) of
Theorem~\ref{uni3} includes the case $k(t)=\sqrt{C_0 t^{\gamma}}$
(when $L_\sharp=0$) provided that $f\in {\mathcal F}_{\rho0,\tau}$
with $\ell^\star\not=0$. Relations (\ref{chi1}) and (\ref{chi3})
show how dramatically changes the two-term asymptotic expansion of
$u_a$ from the result in \cite{gls}. Our approach is completely
different from that in \cite{be,bm,gls,3}, as we use essentially
Karamata's theory.
\end{remark}

We point out that the asymptotic general results stated in the
above theorems do not concern the difference or the quotient of
$u(x)$ and $\psi (d(x))$, as established in \cite{bm}, \cite{bie},
\cite{3}, \cite{rad} for $a=0$ and $b=1$, where $\psi$ is a large
solution of \[\psi ''(r)=f(\psi (r))\qquad\mbox{on}\
(0,\infty)\,.\] For instance, Bieberbach \cite{bie} and Rademacher
\cite{rad} proved that $|u(x)-\psi (d(x))|$ is bounded in a
neighborhood of the boundary. Their result was improved by Bandle
and Ess\'en \cite{be} who showed that $\lim_{d(x)\ri 0}
\left(u(x)-\psi (d(x))\right)=0$.

\medskip
The rest of the paper is organized as follows. In
section~\ref{first}  we collect the notions and properties of
regularly varying functions that are invoked in our proofs. In
section~\ref{second} we prove some auxiliary results including
Lemmas 1 and 2 in \cite{cr4}, which have only been stated there.
In Section~\ref{gyr} we characterize the class ${\mathcal K}$ as
well as its subclasses ${\mathcal K}_{0,\zeta}$ and ${\mathcal
K}_{(01],\tau}$ that appear in Theorems~\ref{uni2} and \ref{uni3}.
Sections~\ref{nup} and \ref{five} are dedicated to the proof of
Theorems~\ref{uni1} and \ref{uni3}.

\section{Preliminaries} \label{sec2}

\subsection{Properties of regularly varying function} \label{first}
The theory of regular variation was instituted in 1930 by Karamata
\cite{kar,kar1} and subsequently developed by himself and many
others. Although Karamata originally introduced his theory in
order to use it in Tauberian theorems, regularly varying functions
have been later applied in several branches of Analysis: Abelian
theorems (asymptotic of series and integrals---Fourier ones in
particular), analytic (entire) functions, analytic number theory,
etc. The great potential of regular variation for probability
theory and its applications was realised by Feller \cite{fe} and
also stimulated by de Haan \cite{ha}. The first monograph on
regularly varying functions was written by Seneta \cite{se}, while
the theory and various applications of the subject are presented
in the comprehensive treatise of Bingham, Goldie and Teugels
\cite{bgt}.

We give here a brief account of the definitions and properties of
regularly varying functions involved in our paper (see \cite{bgt}
or \cite{se} for details).

\begin{definition}\label{def1}
A positive measurable function $Z$ defined on $[A, \infty)$, for
some $A>0$, is called \emph{regularly varying (at infinity) with
index} $q\in \RR$, written $Z\in RV_q$, provided that \[
\lim_{u\to \infty}\frac{Z(\xi u)}{Z(u)}=\xi^q,\qquad \mbox{for
all}\ \xi>0.\] When the index $q$ is zero, we say that the
function is \emph{slowly varying}.
\end{definition}

\begin{remark}\label{tr} Let $Z:[A,\infty)\to (0,\infty)$ be a
measurable function. Then
\begin{enumerate}
\item $Z$ is regularly varying if and only if $\lim_{u\to
\infty}Z(\xi u)/Z(u)$ is finite and positive for each $\xi$ in a
set $S\subset (0,\infty)$ of positive measure (see \cite[Lemma 1.6
and Theorem 1.3]{se}). \item The transformation $Z(u)=u^q L(u)$
reduces regular variation to slow variation. Indeed, $\lim_{u\to
\infty}Z(\xi u)/Z(u)=u^q$ if and only if $\lim_{u\to \infty} L(\xi
u)/L(u)=1$, for every $\xi>0$.
\end{enumerate}
\end{remark}

\begin{example}
Any measurable function on $[A,\infty)$ which has a positive limit
at infinity is slowly varying. The logarithm $\log u$, its
iterates $\log\log u$ ($=\log_2 u$), $\log_m u$ ($=\log\log_{m-1}u
$) and powers of $\log_m u$ are non-trivial examples of slowly
varying functions. Non-logarithmic examples are given by $ {\rm
exp}\left\{(\log u)^{\alpha_1}\right\}$, where $\alpha_1\in (0,1)$
and ${\rm exp} \left\{(\log u)/\log\log u \right\}$.
\end{example}

In what follows $L$ denotes a slowly varying function defined on
$[A,\infty)$. For details on Propositions~\ref{pn}--\ref{kam}, we
refer to \cite{bgt}.

\begin{proposition}[Uniform Convergence Theorem]
\label{pn} The convergence $\frac{L(\xi u)}{L(u)}\to 1$ as $u\to
\infty$ holds uniformly on each compact $\xi$-set in $(0,\infty)$.
\end{proposition}

\begin{proposition}[Representation Theorem]
\label{p1} The function $L(u)$ is slowly varying if and only if it
can be written in the form
\neweq L(u)=M(u){\rm
exp}\left\{\int_B^u\frac{y(t)}{t}\,dt\right\} \quad (u\geq B)
\label{repres}
\endeq
for some $B>A$, where $y\in C[B,\infty)$ satisfies
$\lim_{u\to\infty}y(u)=0$ and $M(u)$ is measurable on $[B,\infty)$
such that $\lim_{u\to \infty}M(u):=\overline M\in (0,\infty)$.
\end{proposition}

The Karamata representation (\ref{repres}) is non-unique because
we can adjust one of $M(u)$, $y(u)$ and modify properly the other
one. Thus, the function $y$ may be assumed arbitrarily smooth, but
the smothness properties of $M(u)$ can ultimately reach those of
$L(u)$. If $M(u)$ is replaced by its limit at infinity $\overline
M>0$, we obtain a slowly varying function $L_0\in C^1[B,\infty)$
of the form \[ L_0(u)=\overline M {\rm exp}\left\{\int_B^u
\frac{y(t)}{t}\,dt\right\} \quad (u\geq B), \] where $y\in
C[B,\infty)$ vanishes at infinity. Such a function $L_0(u)$ is
called a \emph{normalised} slowly varying function.

As an important subclass of $RV_q$, we distinguish $NRV_q$ defined
as
\neweq NRV_q=\left\{
Z\in RV_q:\ \frac{Z(u)}{u^{q}}\ \mbox{is a \emph{normalised}
slowly varying function}\right\}. \label{iso} \endeq

Notice that $L(u)$ given by (\ref{repres}) is asymptotic
equivalent to $L_0(u)$, which has much enhanced properties. For
instance, we see that $ y(u)=\frac{uL_0'(u)}{L_0(u)}$, for all
$u\geq B$. Conversely, any function $L_0\in C^1[B,\infty)$ which
is positive and satisfies
\neweq \lim_{u\to \infty}
\frac{uL_0'(u)}{L_0(u)}=0 \label{tal}
\endeq
is a normalised slowly varying. More generally, if the right hand
side of (\ref{tal}) is $q\in \RR$, then $L_0\in NRV_q$.

\begin{proposition}[Elementary properties
of slowly varying functions] \label{p2} If $L$ is slowly varying,
then
\begin{enumerate}
\item For any $\alpha>0$, $u^\alpha L(u)\to \infty$,
$u^{-\alpha}L(u)\to 0$ as $u\to \infty $; \item $(L(u))^\alpha$
varies slowly for every $\alpha\in \RR$; \item If $L_1$ varies
slowly, so do $L(u)L_1(u)$ and $L(u)+L_1(u)$.
\end{enumerate}
\end{proposition}

From Proposition~\ref{p2} (i) and Remark~\ref{tr} (ii),
$\lim_{u\to \infty}Z(u)=\infty$ (resp., $0$) for any function
$Z\in RV_q$ with $q>0$ (resp., $q<0$).

\begin{remark} \label{mon} Note that the behavior at infinity for a slowly varying
function cannot be predicted. For instance, $$L(u)={\rm
exp}\left\{(\log u)^{1/3} \cos((\log u)^{1/3})\right\}$$ exhibits
infinite oscillation in the sense that $$\liminf_{u\to \infty}
L(u)=0\ \ \mbox{and} \ \ \limsup_{u\to \infty}L(u)=\infty.$$
\end{remark}

\begin{proposition}[Karamata's Theorem; direct half]
\label{p3} Let $Z\in RV_q$ be locally bounded in $[A,\infty)$.
Then
\begin{enumerate}
\item for any $j\geq -(q+1)$,
\neweq \lim_{u\to \infty}\frac{u^{j+1}Z(u)}{\int_A^u x^jZ(x)\,dx}=
j+q+1. \label{da}
\endeq
\item for any $j<-(q+1)$ (and for $j=-(q+1)$ if $\int^\infty
x^{-(q+1)}Z(x)\,dx<\infty$)
\neweq \lim_{u\to \infty}\frac{u^{j+1}Z(u)}{\int_u^\infty x^jZ(x)\,dx}=
-(j+q+1). \label{dad}
\endeq
\end{enumerate}
\end{proposition}

\begin{proposition}[Karamata's Theorem; converse half]
\label{kam} Let $Z$ be positive and locally integrable in
$[A,\infty)$.
\begin{enumerate}
\item If $(\ref{da})$ holds for some $j>-(q+1)$, then $Z\in RV_q$.
\item If $(\ref{dad})$ is satisfied for some $j<-(q+1)$, then
$Z\in RV_q$.
\end{enumerate}
\end{proposition}

For a non-decreasing function $H$ on $\RR$, we define the (left
continuous) inverse of $H$ by \[ H^\leftarrow (y)=\inf\{s:\
H(s)\geq y\}.\]

\begin{proposition}[see Proposition 0.8 in \cite{res}] \label{p08} We have
\begin{enumerate}
\item If $Z\in RV_q$, then $\lim_{u\to \infty}\log Z(u)/\log u=q$.
\item If $Z_1\in RV_{q_1}$ and $Z_2\in RV_{q_2}$ with $\lim_{u\to
\infty}Z_2(u)=\infty$, then \[ Z_1\circ Z_2\in RV_{q_1 q_2}.\]
\item Suppose $Z$ is non-decreasing, $Z(\infty)=\infty$, and $Z\in
RV_q$, $0< q<\infty$. Then
\[ Z^\leftarrow \in RV_{1/q}.\]
\item Suppose $Z_1$, $Z_2$ are non-decreasing and $q$-varying,
$0<q<\infty$. Then for $c\in (0,\infty)$
\[\lim_{u\to \infty} \frac{Z_1(u)}{Z_2(u)}=c\ \ \mbox{if and only if  }
\lim_{u\to \infty}
\frac{Z_1^\leftarrow(u)}{Z_2^\leftarrow(u)}=c^{-1/q}.\]
\end{enumerate}
\end{proposition}

\subsection{Auxiliary results} \label{second}
Based on regular variation theory, we prove here two results that
have only been stated in \cite{cr4}.

\begin{remark}\label{cori1}
If $f\in RV_{\rho+1}$ ($\rho>0$) is continuous, then
\begin{equation} \label{aga} \Xi(u):=\frac{\sqrt{F(u)}}{f(u)\int_u^\infty
[F(s)]^{-1/2}\,ds}\to \frac{\rho}{2(\rho+2)}\ \ \mbox{as } u\to
\infty ,
\end{equation} where $F$ stands for an antiderivative of $f$.
Indeed, by Proposition~\ref{p3}, we have \begin{equation}
\label{jar} \lim_{u\to \infty}
\frac{F(u)}{uf(u)}=\frac{1}{\rho+2}\ \ \mbox{and} \ \ \lim_{u\to
\infty} \frac{u[F(u)]^{-1/2}}{\int_u^\infty
[F(s)]^{-1/2}\,ds}=\frac{\rho}{2}.\end{equation}
\end{remark}

\begin{lemma}[Properties of $h$] \label{aux} If $f\in RV_{\rho+1}$ $(\rho>0)$
is continuous and $k\in {\mathcal   K}$, then $h$ defined by
{\rm \eq{defi}} is a $C^2$-function satisfying 
the following:
\begin{enumerate}
\item[{\rm (i)}] $\di\lim\limits_{t\searrow
0}\frac{h''(t)}{k^2(t)f(h(t)\xi)}=
\frac{2+\rho\ell_1}{\xi^{\rho+1}(2+\rho)}$, for each $\xi>0$;
\item[{\rm (ii)}] $\di\lim\limits_{t\searrow
0}\frac{h(t)h''(t)}{[h'(t)]^2}=\frac{2+\rho \ell_1}{2}$ and
$\di\lim\limits_{t\searrow 0}\frac{\ln k(t)}{\ln
h(t)}=\frac{\rho(\ell_1-1)}{2}$; \item[{\rm (iii)}]
$\di\lim\limits_{t\searrow 0}\frac{h'(t)}{th''(t)}=-\frac{\rho
\ell_1}{2+\rho\ell_1}$ and 
$\di\lim\limits_{t\searrow 0}\frac{h(t)}{t^2 h''(t)}=\frac{\rho^2
\ell_1^2}{2(2+\rho \ell_1)}$; \item[{\rm (iv)}]
$\di\lim\limits_{t\searrow
0}\frac{h(t)}{th'(t)}=\lim\limits_{t\searrow 0}\frac{\ln t}{\ln
h(t)}=-\frac{\rho\ell_1}{2}$; \item[{\rm (v)}]
$\di\lim\limits_{t\searrow 0}t^j h(t)=\infty,$ for all $j>0$,
provided that $k\in {\mathcal   K}_0$. If, in addition, $k\in
{\mathcal K}_{0,\zeta}$ then \(\di \lim_{t\searrow 0}
\frac{1}{-\zeta t^\zeta \ln h(t)}= \lim_{t\searrow
0}\frac{h'(t)}{t^{\zeta+1}h''(t)}=
\frac{-\rho L_\star} {2(\zeta+1)}\,. 
\)
\end{enumerate}
\end{lemma}
\begin{proof} By \eq{defi}, the function $h\in C^2(0,\nu)$, for some $\nu>0$, and
$\lim_{t\searrow 0}h(t)=\infty$.

For any $t\in (0,\nu)$, we have $h'(t)=-k(t)\sqrt{2F(h(t))}$ and
\neweq
h''(t)=k^2(t)f(h(t))\left\{1+2\Xi(h(t)) \left[ \left(
\frac{\int_0^t k(s)\,ds}{k(t)}\right)'-1\right] \right\}.
\label{vox} \endeq Using Remark~\ref{cori1} and $f\in
RV_{\rho+1}$, we reach (i). 

(ii). By (i) and \eq{jar}, we get
\neweq
\lim_{t\searrow 0}\frac{h(t)h''(t)}{[h'(t)]^2}= \lim_{t\searrow
0}\frac{h''(t)}{k^2(t)f(h(t))}\,\frac{h(t)f(h(t))}{2F(h(t))}=
\frac{2+\rho \ell_1}{2},\label{no1} \endeq respectively
\neweq
\lim_{t\searrow 0}\frac{k'(t)}{k(t)}\,\frac{h(t)}{h'(t)}=
\lim_{t\searrow 0}\frac{h(t)f(h(t))}{F(h(t))}\,\frac{-k'(t)\left(
\int_0^t k(s)\,ds\right)}{k^2(t)}\,\Xi(h(t))
=\frac{\rho(\ell_1-1)}{2}.\label{lop}
\endeq

(iii). Using (i) and Remark~\ref{cori1}, we find
\[\lim_{t\searrow 0}\frac{h'(t)}{th''(t)}
=\frac{-2(2+\rho)}{2+\rho\ell_1}\, \lim_{t\searrow
0}\frac{\int_0^t k(s)\,ds}{tk(t)}\,\Xi(h(t))
=\frac{-\rho\ell_1}{2+\rho\ell_1}, \] which, together with
\eq{no1}, implies that
\[
\lim_{t\searrow 0}\frac{h(t)}{t^2 h''(t)}= \lim_{t\searrow 0}
\frac{h(t)h''(t)}{[h'(t)]^2}\,
\left[\frac{h'(t)}{th''(t)}\right]^2=
\frac{\rho^2\ell_1^2}{2(2+\rho\ell_1)}.
\]

(iv). If $\ell_1\not=0$, then by (iii), we have \[ \lim_{t\searrow
0}\frac{h(t)}{th'(t)}= \lim_{t\searrow
0}\frac{h(t)}{t^2h''(t)}\,\frac{th''(t)}{h'(t)}=
\frac{-\rho\ell_1}{2}. \]

If $\ell_1=0$, then we derive
\neweq \lim_{t\searrow 0}\frac{k(t)}{tk'(t)}=
\lim_{t\searrow 0}\frac{k^2(t)}{k'(t)\left(\int_0^t
k(s)\,ds\right)}\, \frac{\int_0^t k(s)\,ds}{tk(t)}=0. \label{mur}
\endeq
This and (\ref{lop}) yield \(\lim_{t\searrow
0}\frac{h(t)}{th'(t)}= 0\), which concludes (iv).

(v). If $k\in {\mathcal   K}_0$, then using (iv), we obtain $
\lim_{t\searrow 0}\ln [t^j h(t)]=\infty$, for all $j>0$.

Suppose $k\in {\mathcal   K}_{0,\zeta}$, for some $\zeta>0$. Then,
$\lim_{t\searrow 0}\frac{\int_0^t k(s)\,ds
}{t^{\zeta+1}k(t)}= \frac{L_\star}{\zeta+1}$ 
and
\neweq
\frac{L_\star}{\zeta+1}= \lim_{t\searrow 0}\frac{\int_0^t
k(s)\,ds}{t^{\zeta+1}k(t)}\, \frac{k^2(t)}{k'(t)\left(\int_0^t
k(s)\,ds\right)}= \lim_{t\searrow
0}\frac{k(t)}{t^{\zeta+1}k'(t)}=\frac{-1}{\zeta}\lim_{t\searrow
0}\frac{1}{t^{\zeta}\ln k(t)}. \label{sol}
\endeq

By \eq{no1}, \eq{lop} and \eq{sol}, we deduce
\[ \lim_{t\searrow 0}\frac{h'(t)}{t^{\zeta+1}h''(t)}=
\lim_{t\searrow 0} \frac{h(t)}{h'(t)t^{\zeta+1}}= \lim_{t\searrow
0}\frac{k'(t)h(t)}{k(t)h'(t)}\,\frac{k(t)}{t^{\zeta+1}k'(t)} =
\frac{-\rho L_\star}{2(\zeta+1)}.\] This completes the proof of
the lemma.
\end{proof}

Let $\tau>0$ be arbitrary and $f$ be as in Remark~\ref{cori1}. For
$u>0$ sufficiently large, we define
\neweq T_{1,\tau}(u)=\left[\frac{\rho}{2(\rho+2)}-\Xi(u)\right](\ln
u)^\tau\ \ \mbox{and} \ \ T_{2,\tau}(u)=\left[\frac{f(\xi_0
u)}{\xi_0 f(u)}-\xi_0^\rho\right] (\ln u)^\tau. \label{yz}\endeq

\begin{remark}\label{obs} When $f(u)=C u^{\rho+1}$,
we have $T_{1,\tau}(u)=T_{2,\tau}(u)=0$.
\end{remark}

\begin{lemma}\label{xua} Assume that $f\in {\mathcal  F}_{\rho\eta}$ (where $-\rho-2<\eta\leq 0$).
The following hold:
\begin{enumerate}
\item[(i)] If $f\in {\mathcal  F}_{\rho0,\tau}$, then
\[\lim_{u\to\infty}T_{1,\tau}(u)=\frac{-\ell^\star}{(\rho+2)^2}\ \
\mbox{and}\ \ \lim_{u\to \infty} T_{2,\tau}(u)=\xi_0^\rho
\ell^\star \ln \xi_0.\] \item[(ii)] If $f\in {\mathcal
F}_{\rho\eta}$ with $\eta\not=0$, then $$\lim_{u\to
\infty}T_{1,\tau}(u)=\lim_{u\to \infty}T_{2,\tau}(u)=0 .$$
\end{enumerate}
\end{lemma}
\begin{proof}
Using the second limit in (\ref{jar}), we obtain
\[\lim_{u\to \infty}T_{1,\tau}(u)= \frac{\rho}{2}\lim_{u\to
\infty}\frac{\frac{\rho}{2(\rho+2)}\int_u^\infty
[F(s)]^{-1/2}\,ds-\frac{\sqrt{F(u)}}{f(u)}} {u [F(u)]^{-1/2}\,(\ln
u)^{-\tau}}. \] By L'Hospital's rule, we arrive at
\[\lim_{u\to \infty}T_{1,\tau}(u)=\lim_{u\to \infty} \left[\frac{\rho+1}{\rho+2}-
\frac{F(u)f'(u)}{f^2(u)} \right](\ln u)^\tau :=\lim_{u\to
\infty}Q_{1,\tau}(u).
\]
A simple calculation shows that, for $u>0$ large,
\[\begin{aligned} Q_{1,\tau}(u)&=
\frac{(\ln u)^\tau}{\rho+2}\left[\rho+1-\frac{uf'(u)}{f(u)}
\right]+ \frac{uf'(u)}{f(u)}\left[\frac{1}{\rho+2}-
\frac{F(u)}{uf(u)}\right](\ln u)^\tau\\
&=:\frac{1}{\rho+2}\,Q_{2,\tau}(u)
+\frac{uf'(u)}{f(u)}\,Q_{3,\tau}(u).
\end{aligned}\]
Since (\ref{form}) holds with $\phi\in RV_\eta$ or $-\phi\in
RV_\eta$, we can assume $B>0$ such that $\phi\not= 0$ on
$[B,\infty)$. For any $u>B$, we have $ Q_{2,\tau}(u)= -\phi(u)(\ln
u)^\tau $ and \[ Q_{3,\tau}(u)= \widetilde C\, \frac{(\ln
u)^\tau}{uf(u)}+\frac{\int_B^u f(s)\phi(s)\,ds}{(\rho+2)
uf(u)\phi(u)}\,\phi(u)(\ln u)^\tau,\] where $\widetilde C\in \RR$
is a constant. Since either $f\phi\in RV_{\rho+\eta+1}$ or
$-f\phi\in RV_{\rho+\eta+1}$, by Proposition~\ref{p3}, \[
\lim_{u\to \infty}\frac{uf(u)\phi(u)}
{\int_B^uf(x)\phi(x)\,dx}=\rho+\eta+2.\] If (i) holds, then
$\lim_{u\to \infty}Q_{2,\tau}(u)=-\ell^\star$ and $\lim_{u\to
\infty}Q_{3,\tau}(u)=\ell^\star(\rho+2)^{-2}$. Thus,
\[ \lim_{u\to \infty} T_{1,\tau}(u)= \lim_{u\to
\infty}Q_{1,\tau}(u)=-\ell^\star/(\rho+2)^2 .\] If (ii) holds,
then by Proposition~\ref{p2}, we have $\lim_{u\to \infty}(\ln
u)^\tau \phi(u)=0$. It follows that \[ \lim_{u\to
\infty}Q_{2,\tau}(u)=\lim_{u\to \infty}Q_{3,\tau}(u)=0
\] which yields $\lim_{u\to \infty}T_{1,\tau}(u)=0$. Note that the
proof is finished if $\xi_0=1$, since $T_{2,\tau}(u)=0$ for each
$u>0$.

Arguing by contradiction, let us suppose that $\xi_0\not=1$. Then,
by (\ref{form}), \[ T_{2,\tau}(u)= \xi_0^\rho\left[{\rm
exp}\left\{\int_u^{\xi_0 u}\frac{\phi(t)}{t}\,dt\right\}-1\right]
(\ln u)^\tau,\quad \forall u>B/\xi_0. \] But, $\lim_{u\to
\infty}\phi(us)/s=0$, uniformly with respect to $s\in [\xi_0,1]$.
So \[ \lim_{u\to \infty}\int_u^{\xi_0 u}\frac{\phi(t)}{t}\,dt=
\lim_{u\to \infty} \int_1^{\xi_0}\frac{\phi(su)}{s}\,ds =0 \]
which leads to \[ \lim_{u\to
\infty}T_{2,\tau}(u)=\xi_0^\rho\lim_{u\to \infty}
\left(\int_u^{\xi_0 u}\frac{\phi(t)}{t}\,dt\right)(\ln u)^\tau. \]

If (i) occurs, then by Proposition~\ref{pn}, we have \[
 \lim_{u\to \infty}T_{2,\tau}(u)=\xi_0^\rho
\lim_{u\to \infty}(\ln u)^\tau \phi(u) \int_1^{\xi_0}
\frac{\phi(tu)}{\phi(u)}\,\frac{dt}{t} =\xi_0^\rho\ell^\star\ln
\xi_0.\]

If (ii) occurs, then by Proposition~\ref{p2}, we infer that
\[ \lim_{u\to \infty}T_{2,\tau}(u)=
\frac{-\xi_0^\rho}{\tau}\lim_{u\to \infty} \left[\phi(\xi_0
u)-\phi(u)\right] (\ln u)^{\tau+1}=0.\] The proof of
Lemma~\ref{xua} is now complete.  \end{proof}

\begin{lemma} \label{goal}
If $k\in {\mathcal K}_{(01],\tau}$ and $f$ satisfies either {\rm
(i)} or {\rm (ii)} of Theorem~{\rm \ref{uni3}}, then
\neweq \label{pik} {\mathcal H}(t):=(-\ln t)^\tau\left(1-\frac{k^2(t)f(\xi_0 h(t))}
{\xi_0 h''(t)}\right)\to \rho \widetilde \chi\ \ \mbox{as
}t\searrow 0,
\endeq
where $\widetilde \chi$ is defined by {\rm \eq{hii}}.
\end{lemma}

\begin{proof} Using (\ref{vox}), we write $ {\mathcal
H}(t)=\frac{k^2(t)f(h(t))}{h''(t)}\sum_{i=1}^3 {\mathcal
H}_{i}(t)$, for $t>0$ small, where
\[\left\{\begin{aligned}
{\mathcal H}_{1}(t)&:= 2\Xi(h(t))(-\ln t)^\tau \left[
\left(\frac{\int_0^t k(s)\,ds}{k(t)}\right)'-\ell_1\right],\\
{\mathcal H}_{2}(t)&:= 2(1-\ell_1)\left(\frac{-\ln t}{\ln
h(t)}\right)^\tau T_{1,\tau}(h(t))\ \ \mbox{and}\ \ {\mathcal
H}_{3}(t):= - \left(\frac{-\ln t}{\ln h(t)}\right)^\tau
T_{2,\tau}(h(t)).\end{aligned}\right.\]

By Remark~\ref{cori1}, we find \(\lim_{t\searrow 0} {\mathcal
H}_{1}(t)=\rho L_\sharp/(\rho+2)\).

\emph{Case} (i) (that is, $f\in {\mathcal  F}_{\rho\eta}$ with
$\eta L_\sharp\not=0$). By Lemmas~\ref{aux} and \ref{xua}, it
turns out that \[ \lim_{t\searrow 0}{\mathcal
H}_{2}(t)=\lim_{t\searrow 0}{\mathcal H}_{3}(t)=0 \ \ \mbox{and} \
\ \lim_{t\searrow 0}{\mathcal H}(t)=\frac{\rho
L_\sharp}{2+\rho\ell_1}=:\rho \widetilde \chi.
\]

\emph{Case} (ii) (that is, $f\in {\mathcal F}_{\rho0,\tau}$ with
$[\ell^\star(\ell_1-1)]^2+L_\sharp^2\not=0$). By Lemmas~\ref{aux}
and \ref{xua}, we get
\[ \lim_{t\searrow 0} {\mathcal H}_{2}(t)=
\frac{-2(1-\ell_1)\ell^\star}{(\rho+2)^2}\, \left(\frac{\rho
\ell_1}{2}\right)^\tau \ \ \mbox{and}\ \ \lim_{t\searrow
0}{\mathcal H}_{3}(t)= \frac{-\ell^\star(2+\rho\ell_1)}{(2+\rho)}
\left(\frac{\rho\ell_1}{2}\right)^\tau \ln \xi_0.\] Thus, we
arrive at
\[ \lim_{t\searrow 0}{\mathcal H}(t)= \frac{\rho
L_\sharp}{2+\rho\ell_1}-\ell^\star
\left(\frac{\rho\ell_1}{2}\right)^\tau
\left[\frac{2(1-\ell_1)}{(\rho+2)(2+\rho\ell_1)}+\ln
\xi_0\right]=:\rho \widetilde \chi.\] This finishes the proof.
\end{proof}

\section{Characterization of ${\mathcal  K}$ and its subclasses}\label{gyr}

Definition~\ref{def1} extends to \emph{regular variation at the
origin}. We say that $Z$ is regularly varying (on the right) at
the origin with index $q$ (and write, $Z\in RV_q(0+))$ if
$Z(1/u)\in RV_{-q}$. Moreover, by $Z\in NRV_q(0+)$ we mean that
$Z(1/u)\in NRV_{-q}$. The meaning of $NRV_q$ is given by
(\ref{iso}).

\begin{proposition} \label{lg1} We have $k\in {\mathcal K}_{(01]}$ if and only
if $k$ is non-decreasing near the origin and $k$ belongs to
$NRV_{\alpha}(0+)$ for some $\alpha\geq 0$ (where
$\alpha=1/\ell_1-1$).
\end{proposition}
\begin{proof}
If $k\in {\mathcal K}_{(01]}$, then from the definition
\[\lim_{t\to 0^+}\frac{\int_0^tk(s)ds}{k(t)}\Big/t= \lim_{t\to
0^+}\Big(\frac{\int_0^tk(s)ds}{k(t)}\Big)'=\ell_1,\] which implies
that
\[ \lim_{u\to \infty}\frac{u\,\frac{d}{du}k(1/u)}{k(1/u)}=\lim_{t\to 0^+}
\frac{-tk'(t)}{k(t)}=\frac{\ell_1-1}{\ell_1}.
\] Thus $k(1/u)$ belongs to $NRV_{1-1/\ell_1}$.
Conversely, if $k$ belongs to $NRV_\alpha(0+)$ with $\alpha\geq
0$, then $k$ is a positive $C^1$-function on some interval
$(0,\nu)$ and
\begin{equation} \label{opt}\lim_{t\to 0^+}
\frac{tk'(t)}{k(t)}=\alpha.
\end{equation}
By Proposition~\ref{p3}, we deduce \neweq \label{opt2} \lim_{t\to
0^+}\frac{\int_0^t k(s)\,ds}{t k(t)}=
\lim_{u\to\infty}\frac{\int_u^\infty
x^{-2}k(1/x)dx}{u^{-1}k(1/u)}=\frac{1}{1+\alpha}.\endeq Combining
\eq{opt} and \eq{opt2}, we get $\lim_{t\to 0^+}\left(\int_0^t
k(s)\,ds/k(t)\right)'=1/(1+\alpha)$. If, in addition, $k$ is
non-decreasing near $0$, then $k\in {\mathcal K}$ with
$\ell_1=1/(1+\alpha)$. Note that by \eq{opt}, $k$ is increasing
near the origin if $\alpha>0$; however, when $k$ is slowly varying
at $0$, then we cannot draw any conclusion about the monotonicity
of $k$ near the origin (see Remark~\ref{mon}).
\end{proof}

\begin{remark} By Propositions~\ref{lg1} and \ref{pn}, we deduce $k\in
{\mathcal K}_{(01]}$ if and only $k$ is of the form
\neweq \label{lif}
k(t)=c_0 t^\alpha
\exp\left\{\int_t^{c_1}\frac{E(y)}{y}\,dy\right\} \ \ (0<t<c_1), \
\mbox{for some } 0\leq \alpha(=1/\ell_1-1)
\endeq
where $c_0,c_1>0$ are constants, $E\in C[0,c_1)$ with $E(0)=0$ and
(only for $\ell_1=1$) $E(t)\leq \alpha$.
\end{remark}

\begin{proposition} \label{lg2} We have $k\in {\mathcal K}_{(01],\tau}$ if
and only if $k$ is of the form {\rm \eq{lif}} where, in addition,
\neweq
\label{lif2} \lim_{t\searrow 0} (-\ln t)^\tau E(t)=\ell_\sharp\in
\RR \ \ \mbox{with } \ell_\sharp=(1+\alpha)^2 L_\sharp.
\endeq
\end{proposition}
\begin{proof}
Suppose $k$ satisfies \eq{lif} and \eq{lif2}. A simple calculation
leads to
\neweq
\lim_{t\searrow 0}(-\ln t)^\tau \left[\frac{1-\ell_1}{\ell_1}-
\frac{tk'(t)}{k(t)}\right]=\lim_{t\searrow 0}(-\ln t)^\tau
E(t)=\ell_\sharp. \label{st1}
\endeq By L'Hospital's rule, we find \begin{equation}
\begin{aligned}
\lim_{t\searrow 0}(-\ln t)^\tau \left[\ell_1-\frac{\int_0^t
k(s)\,ds}{tk(t)} \right] &= \lim_{t\searrow 0}\frac{(\ell_1-1)
+\ell_1 tk'(t)/k(t)}{(-\ln
t)^{-\tau}\left[1+\frac{tk'(t)}{k(t)}-\frac{\tau} {\ln
t}\right]}\\
& = -\ell_1^2 \lim_{t\searrow 0}(-\ln t)^\tau
\left[\frac{1-\ell_1}{\ell_1} -\frac{tk'(t)}{k(t)}\right]=
\frac{-\ell_\sharp}{(\alpha+1)^2}.
\end{aligned} \label{st2}
\end{equation}
We see that, for each $t>0$ small,
\begin{equation}
\left(\frac{\int_0^t k(s)\,ds} {k(t)}\right)'-\ell_1 =
\frac{tk'(t)}{k(t)} \left[ \ell_1-\frac{\int_0^t
k(s)\,ds}{tk(t)}\right]+\ell_1
\left[\frac{1-\ell_1}{\ell_1}-\frac{tk'(t)}{k(t)}\right].
\label{st3}
\end{equation}
By (\ref{st1})--(\ref{st3}), we infer that $k\in {\mathcal
K}_{(01],\tau}$ with $L_\sharp=\ell_\sharp/(1+\alpha)^2$.

Conversely, if $k\in {\mathcal  K}_{(01],\tau} $, then $k$ is of
the form \eq{lif}. Moreover, we have
\neweq \lim_{t\searrow 0}(-\ln t)^\tau
\left(\frac{\int_0^t k(s)\,ds}{tk(t)}-\ell_1
\right)=\lim_{t\searrow 0}\frac{ \left(\int_0^t k(s)\,ds/k(t)
\right)'-\ell_1}{(-\ln t)^{-\tau}\left(1-\frac{\tau}{\ln
t}\right)}=L_\sharp. \label{st4}
\endeq
By (\ref{st3}) and (\ref{st4}), we deduce \[ L_\sharp=-\alpha
L_\sharp+\frac{1}{\alpha+1}\,\lim_{t\searrow 0} (-\ln t)^\tau
E(t).\] Consequently, $\lim_{t\searrow 0}(-\ln t)^\tau
E(t)=(1+\alpha)^2 L_\sharp$. Hence, \eq{lif2} holds. \end{proof}

\begin{proposition} \label{mat} We have $k\in {\mathcal K}_0$ if and only if $k$ is of
the form
\neweq \label{lig} k(t)=d_0 \left(\exp\left\{-\int_t^{d_1} \frac{dx}{x {\mathcal W}(x)}\right\}
\right)' \ \ (0<t<d_1), \endeq where $d_0,d_1>0$ are constants and
$0<{\mathcal W} \in C^1(0,d_1)$ satisfies $\lim_{t\searrow 0}
{\mathcal W}(t)=\lim_{t\searrow 0}t\,{\mathcal W}'(t)=0$.
\end{proposition}

\begin{proof} If $k\in {\mathcal K}_0$, then we set
\neweq \label{defw} {\mathcal W}(t)=\frac{\int_0^t k(s)\,ds}{t k(t)},
\ \ \mbox{for } t \in (0, d_1).\
\endeq
Hence, $\lim_{t\searrow 0} {\mathcal W}(t)=0$ and, for $t>0$
small,
\[ t {\mathcal W}'(t)=\left(\frac{\int_0^t
k(s)\,ds}{k(t)}\right)'-\frac{\int_0^t k(s)\,ds}{tk(t)}.
\]
It follows that $\lim_{t\searrow 0} t{\mathcal W}'(t)=0$. By
\eq{defw}, we find
\[ \int_t^{d_1} \frac{dx}{x {\mathcal W(x)}}=\ln \left(\int_0^{d_1}
k(s)\,ds \right)-\ln \left(\int_0^t k(s)\,ds \right),\ \ t\in
(0,d_1)
\]
so that \eq{lig} is fulfilled. Conversely, if \eq{lig} holds, then
$\lim_{t\to 0} \int_t^{d_1} \frac{dx}{x {\mathcal W}(x)}=\infty$
and \neweq \label{lig3} \int_0^t k(s)\,ds=d_0
\exp\left\{-\int_t^{d_1}\frac{dx}{x {\mathcal W}(x)}\right\}=tk(t)
{\mathcal W}(t),\ \ t\in (0,d_1). \endeq This, together with the
properties of ${\mathcal W}$, shows that $k\in {\mathcal K}_0$.
\end{proof}

\begin{proposition} \label{lok} We have $k\in {\mathcal K}_{0,\zeta}$ if and only if $k$ is
of the form {\rm \eq{lig}} where, in addition,
\neweq \label{lig2} \lim_{t\searrow 0} t^{1-\zeta} {\mathcal W}'(t)=-\ell_\star \ \
\mbox{with } -\ell_\star=\zeta L_\star/(1+\zeta). \endeq
\end{proposition}
\begin{proof} If $k\in {\mathcal K}_{0,\zeta} $, then \eq{lig} and
\eq{lig3} are fulfilled. Therefore,
\[L_\star=\lim_{t\searrow 0}\frac{(t{\mathcal W}(t))'}{t^\zeta}=\lim_{t\searrow 0}\frac{{\mathcal 
W}(t)+
t\,{\mathcal W}'(t)}{t^\zeta} \ \ \mbox{and  }
\frac{L_\star}{\zeta+1} =\lim_{t\searrow 0} \frac{\int_0^{t}
k(s)\,ds}{k(t)t^{\zeta+1}}=\lim_{t\searrow 0} \frac{{\mathcal
W}(t)}{t^{\zeta}},\] from which \eq{lig2} follows. Conversely, if
\eq{lig} and \eq{lig2} hold, then $\lim_{t\searrow 0} {\mathcal
W}(t)/t^\zeta=-\ell_\star/\zeta$. By \eq{lig3}, we infer that
\[ \frac{1}{t^\zeta}\,\left(\frac{\int_0^t
k(s)\,ds}{k(t)}\right)'=\frac{1}{t^\zeta} ({\mathcal W}(t)+t
{\mathcal W}'(t))\to \frac{-\ell_\star(\zeta+1)}{\zeta} \ \
\mbox{as } t\searrow 0.
\]
Thus, $k\in {\mathcal K}_{0,\zeta} $ with
$L_\star=-\ell_\star(\zeta+1)/\zeta$.
\end{proof}

\begin{remark}\label{tur} If $k\in {\mathcal K}_0$ or
$k\in {\mathcal K}_{(01],\tau}$ with
$(1-\ell_1)^2+L_\sharp^2\not=0$, then
\neweq \lim_{t\searrow 0}\frac{k'(t)}{k(t)t^{\theta-1}}=\infty,\quad
\mbox{for every}\ \theta>0. \label{um} \endeq Indeed, if $k\in
{\mathcal K}_0$, then $\lim_{t\searrow
0}\frac{tk'(t)}{k(t)}=\infty$. Assuming that $k\in {\mathcal
K}_{(01],\tau}$, we deduce \eq{um} from \eq{opt} when
$\ell_1\not=1$, otherwise from \eq{lif2} when $L_\sharp\not=0$
since
\[ \lim_{t\searrow
0}\frac{k'(t)}{k(t)t^{\theta-1}}=\lim_{t\searrow 0}
-E(t)t^{-\theta}= -L_\sharp \,\lim_{t\searrow 0}
\frac{t^{-\theta}}{(-\ln t)^\tau}=\infty. \]
\end{remark}

\begin{definition}[see \cite{res}] \label{deli}
A non-decreasing function $U$ is $\Gamma$-varying at $\infty$ if
$U$ is defined on an interval $(A,\infty)$, $\lim_{x\to
\infty}U(x)=\infty$ and there is $g:(A,\infty)\to (0,\infty)$ such
that
\[ \lim_{y\to \infty}\frac{U(y+\lambda g(y))}{U(y)}=e^\lambda,\
\forall \lambda\in \RR .\]
\end{definition}
The function $g$ is called an {\em auxiliary function} and is
unique up to asymptotic equivalence.

\begin{remark} \label{see} Under the assumptions of Theorem~{\rm \ref{uni1}}, we
have
\begin{enumerate}
\item[{\rm (a)}] Suppose $\lim_{t\searrow 0} \left(\int_0^t
k(s)\,ds\right)^2 r(t)=1$ and let $\widehat f(u)$ be such that
$\lim_{u\to \infty} \widehat f(u)/f(u)=1$ and
$j(u):=\widehat{f}(u)/u$ is non-decreasing for $u>0$ large. Then
$\lim_{t\searrow 0}\widehat{\ph}(t)/\ph(t)=1$, where $\ph(t)$ is
given by {\rm \eq{defp}} and $\widehat{\ph}(t)=j^\leftarrow(r(t))$
for $t>0$ small. \item[{\rm (b)}] If $k\in {\mathcal K}$ with
$\ell_1\not=0$, then $\ph(1/u)\in RV_{2/(\rho\ell_1)}$. \item[{\rm
(c)}] If $k\in {\mathcal K}_0$, then $\ph(1/u)$ is
$\Gamma$-varying at $u=\infty$ with auxiliary function
$\di\frac{\rho u^2 \int_0^{1/u}k(s)\,ds}{2k(1/u)}$.
\item[{\rm (d)}] $\lim_{t\searrow 0}
\ph(t)/h(t)=[2(\rho+2)/\rho^2]^{-1/\rho}$, where $h(t)$ is given
by {\rm \eq{defi}}.
\end{enumerate}
Indeed, by Proposition~\ref{p08} we find $(f(u)/u)^\leftarrow \in
RV_{1/\rho}$ and $\lim_{u\to
\infty}(f(u)/u)^\leftarrow/j^\leftarrow(u) =1$. Then, by
Proposition~\ref{pn} we deduce (a). We see that (b) follows by
Proposition~\ref{p08} since $\left(\int_0^{1/u}
k(s)\,ds\right)^{-2}\in RV_{2/\ell_1}$ (cf. Proposition~\ref{lg1})
and $f(u)/u\in RV_\rho$. If $k\in {\mathcal K}_0$, then by
Proposition~\ref{mat} and \cite[p. 106]{res}, we get
$\left(\int_0^{1/u}k(s)\,ds\right)^{-2}$ is $\Gamma$-varying at
$u=\infty$ with auxiliary function $u{\mathcal W}(1/u)/2$. By
\cite[p. 36]{res}, we conclude (c). Notice that
$Y(u):=\left(1/\int_u^\infty [2F(s)]^{-1/2}\,ds\right)^2\in
RV_\rho$ and $Y(h(t))=\left(\int_0^t k(s)\,ds\right)^{-2}$ for
$t>0$ small. We have $\lim_{u\to
\infty}f(u)/[uY(u)]=2(\rho+2)/\rho^2$ (cf. Remark~\ref{cori1}). By
Proposition~\ref{p08}, we achieve (d).
\end{remark}

\section{Proof of Theorem~\ref{uni1}} \label{nup}

Fix $a\in (-\infty,\lambda_{\infty,1})$. By
\cite[Theorem~1.1]{cr2}, equation (\ref{0a}) has at least a large
solution.

In what follows, we will prove that \eq{comp} holds for any large
solution. Hence, a standard argument leads to the uniqueness (see,
for instance, \cite{gls} or \cite{cr2}).

By virtue of Remark~\ref{see} (d), it is enough to demonstrate
\eq{fjk}. Let $u_a$ denote an arbitrary large solution of
(\ref{0a}). Fix $\ep\in (0,1/2)$ and choose $\delta >0$ such that
\begin{enumerate}
\item[(i)] $d(x)$ is a $C^2$ function on the set $\{x\in \Omega:\
d(x)<\delta\}$; \item[(ii)] $ k$ is non-decreasing on
$(0,\delta)$; \item[(iii)] $1-\ep< b(x)/k^2(d(x))<1+\ep$, \quad
$\forall x\in \Omega$ with $0<d(x)<\delta$ (since \eq{B} holds);
\item[(iv)] $h'(t)<0$ and $h''(t)>0$ for each $t\in (0,\delta)$
(cf. Lemma~\ref{aux}).
\end{enumerate}
Define $\xi^\pm=\left[\frac{2+\ell_1\rho}{(1\mp 2\ep)(2+\rho)}
\right]^{1/\rho}$ and $u^\pm (x)=\xi^\pm h(d(x))$, for any $x$
with $d(x)\in (0,\delta)$.

The proof of (\ref{fjk}) will be divided into three steps:

\emph{Step} 1. There exists $\delta_1\in (0,\delta)$ small such
that
\begin{equation}
\left\{\begin{aligned} & \Delta u^++a u^+-(1-\ep)k^2(d)f(u^+)\leq
0, \quad \forall x
\mbox{ with } d(x)\in (0,\delta_1)\\
& \Delta u^-+a u^--(1+\ep)k^2(d)f(u^-) \geq  0, \quad \forall x
\mbox{ with } d(x)\in (0,\delta_1).
\end{aligned} \right.
\label{y1}
\end{equation}

Indeed, for every $x\in \Omega$ with $0<d(x)<\delta$, we have
\begin{equation} \begin{aligned}
& \Delta u^\pm + a u^\pm-(1\mp\ep)k^2(d)f(u^\pm)\\
 & \quad =\xi^\pm h''(d) \left(1+a\frac{h(d)}{h''(d)}+ \Delta
d\,\frac{h'(d)}{h''(d)}-(1\mp\ep)\frac{k^2(d)f(u^\pm)}{\xi^\pm
h''(d)}\right)=:\xi^\pm h''(d)B^\pm(d).
\end{aligned} \label{nil}\end{equation}
By Lemma~\ref{aux}, we deduce $\lim_{d\searrow 0}B^\pm(d)=\mp
\ep/(1\mp 2\ep) $, which proves \eq{y1}.

\emph{Step} 2. There exists $M^+$, $\delta^+>0$ such that \[
u_a(x)\leq u^+(x)+M^+,\quad \forall x\in \Omega\ \mbox{with}\
0<d<\delta^+.
\] For $x\in \Omega$ with $d(x)\in (0,\delta_1)$, we define
$\Psi_x(u)=au-b(x)f(u)$ for each $u>0$. By Lemma~\ref{aux},
\neweq \lim_{d(x)\searrow 0}\frac{b(x)f(u^+(x))}{u^+(x)}=
\lim_{d\searrow 0}\frac{k^2(d)f(u^+)}{\xi^+ h''(d)}
\,\frac{h''(d)}{h(d)}=\infty. \label{emul} \endeq From this and
$(A_1)$, we infer that there exists $\delta_2\in (0,\delta_1)$
such that, for any $x$ with $0<d(x)<\delta_2$, \[u\longmapsto
\Psi_x(u)\ \mbox{is decreasing on some interval } (u_x,\infty) \
\mbox{with } 0<u_x<u^+(x).\] Hence, for each $M>0$, we have
\neweq \Psi_x(u^+(x)+M)\leq \Psi_x(u^+(x)),
\quad \forall x\in \Omega \mbox{ with } 0<d(x)<\delta_2.
\label{z1}\endeq

Fix $\sigma\in (0,\delta_2/4)$ and set $ {\mathcal  N}_{\sigma}:=
\{x\in \Omega:\ \sigma<d(x)<\delta_2/2\}$.

We define $u^*_{\sigma}(x)=u^+(d-\sigma,s)+M^+$, where $(d,s)$ are
the local coordinates of $x\in {\mathcal N}_\sigma$. We choose
$M^+>0 $ large enough such that
\[u^*_{\sigma}(\delta_2/2,s)=u^+(\delta_2/2-\sigma,s)+M^+
\geq u_a(\delta_2/2,s), \quad \forall \sigma\in (0,\delta_2/4)
\mbox{ and } \forall s\in \partial\Omega. \]

By (ii), (iii), (\ref{y1}) and (\ref{z1}), we obtain
\[\begin{aligned}
 -\Delta u^*_\sigma (x) &\geq  a u^+(d-\sigma,s)-(1-\ep)k^2(d-\sigma)f(u^+(d-\sigma,s)) \\
& \geq  a u^+(d-\sigma,s)-b(x)f(u^+(d-\sigma,s))\\
&\geq a (u^+(d-\sigma,s)+M^+)-b(x)f(u^+(d-\sigma,s)+M^+)\\
&= au^*_\sigma(x)-b(x)f(u^*_\sigma(x))\quad \mbox{in } {\mathcal
N}_{\sigma}.
\end{aligned}\]
So, uniformly with respect to $\sigma$, we have
\neweq \Delta u_{\sigma}^*(x)+au_\sigma^*(x)\leq b(x)
f(u_\sigma^*(x))\quad \mbox{in } {\mathcal  N}_{\sigma}.
\label{emu2}
\endeq Since $u_\sigma^*(x)\to \infty$ as $d\searrow \sigma$,
from \cite[Lemma~2.1]{cr2}, we get $ u_a\leq u_\sigma^*$ in
${\mathcal N}_{\sigma}$, for every $\sigma\in (0,\delta_2/4)$.
Letting $\sigma\searrow 0$, we achieve the assertion of Step 2
(with $\delta^+\in (0,\delta_2/2)$ arbitrarily chosen).

\emph{Step} 3. There exists $M^-$, $\delta^->0$ such that
\neweq u_a(x)\geq u^-(x)-M^-,\quad \forall x=(d,s)\in \Omega\quad
\mbox{with}\ 0<d<\delta^-. \label{nb} \endeq

For every $r\in (0,\delta)$, define $\Omega_r= \{x\in \Omega:\
0<d(x)<r\}$.

Fix $\sigma\in (0,\delta_2/4)$. We define $v^*_\sigma(x)= \lambda
u^-(d+\sigma,s) $ for $x=(d,s)\in \Omega_{\delta_2/2}$, where
$\lambda\in (0,1)$ is chosen small enough such that
\neweq
v_{\sigma}^*(\delta_2/4,s)= \lambda u^-(\delta_2/4+\sigma,s) \leq
u_a(\delta_2/4,s), \quad \forall \sigma\in (0,\delta_2/4),\
\forall s\in \partial\Omega. \label{la}\endeq Notice that
$\limsup_{d\searrow 0}(v_\sigma^*-u_a)(x)=-\infty$. By (ii),
(iii), (\ref{y1}) and $(A_1)$, we have
\[\begin{aligned}
 \Delta v^*_\sigma (x)+av_\sigma^*(x)  &=
\lambda(\Delta u^-(d+\sigma,s)+au^-(d+\sigma,s))\\
 &\geq  \lambda(1+\ep) k^2(d+\sigma)f(u^-(d+\sigma,s))
 \geq  (1+\ep) k^2(d)f(\lambda u^-(d+\sigma,s))\\
 &\geq b(x)f(v^*_\sigma(x)),\quad
\forall x=(d,s)\in \Omega_{\delta_2/4}.
\end{aligned}\]
Using \cite[Lemma~2.1]{cr2}, we derive $v_\sigma^*\leq u_a$ in
$\Omega_{\delta_2/4}$. Letting $\sigma\searrow 0$, we get
\neweq
\lambda u^-(x)\leq u_a(x),\quad \forall x\in \Omega_{\delta_2/4}.
\label{em}
\endeq
By Lemma~\ref{aux}, $ \lim_{d\searrow 0} k^2(d)f(\lambda^2 u^-)/
u^-=\infty $. Thus, there exists $\widetilde\delta\in
(0,\delta_2/4)$ such that
\neweq  k^2(d)f(\lambda^2 u^-)/u^-
\geq \lambda^2|a|,\quad \forall x\in \Omega\ \mbox{with}\ 0<d\leq
\widetilde\delta. \label{et}\endeq Choose $\delta_*\in
(0,\widetilde \delta)$, sufficiently close to $\widetilde\delta$,
such that
\neweq h(\delta_*)/h(\widetilde\delta)<1+\lambda.
\label{emu}\endeq

For each $\sigma\in (0,\widetilde\delta-\delta_*)$, we define
\(z_\sigma(x)=u^-(d+\sigma,s)-(1-\lambda)u^-(\delta_*,s)\), where
$x=(d,s)\in \Omega_{\delta_*}$. We prove that $z_\sigma$ is
positive in $\Omega_{\delta_*}$ and
\neweq \Delta z_\sigma+a z_\sigma\geq b(x)f(z_\sigma)
\quad \mbox{in } \Omega_{\delta_*}. \label{hy}\endeq By (iv),
$u^-(x)$ decreases with $d$ when $d<\widetilde\delta$. This and
(\ref{emu}) imply that
\neweq \label{bim} 1+\lambda>\frac{u^-(\delta_*,s)}{u^-(\widetilde{\delta},s)} \geq
\frac{u^-(\delta_*,s)}{u^-(d+\sigma,s)},\quad \forall x=(d,s)\in
\Omega_{\delta_*}.\endeq Hence, \begin{equation} \label{ja}
z_\sigma(x)=u^-(d+\sigma,s)
\left(1-\frac{(1-\lambda)u^-(\delta_*,s)}{u^-(d+\sigma,s)}\right)
\geq \lambda^2 u^-(d+\sigma,s)>0,\ \ \forall x\in
\Omega_{\delta_*}.\end{equation}

By (\ref{y1}), (ii) and (iii), we see that (\ref{hy}) follows if
\begin{equation}
(1+\ep)k^2(d +\sigma)\left[
f(u^-(d+\sigma,s))-f(z_\sigma(d,s))\right]\geq
a(1-\lambda)u^-(\delta_*,s),\quad \forall (d,s)\in
\Omega_{\delta_*}. \label{inew}
\end{equation}
The Lagrange mean value theorem and $(A_1)$ show that \neweq
\label{bem} f(u^-(d+\sigma,s))-f(z_\sigma(d,s)) \geq (1-\lambda)
u^-(\delta_*,s)f(z_\sigma(x))/z_\sigma(x) \endeq which, combined
with \eq{et} and \eq{ja}, proves \eq{inew}.

Notice that $\limsup_{d\searrow 0}(z_\sigma- u_a)(x)=-\infty$. By
(\ref{em}), we have
\[z_\sigma(x)=u^-(\delta_*+\sigma,s)-(1-\lambda)u^-(\delta_*,s) \leq
\lambda u^-(\delta_*,s) \leq u_a(x),\quad \forall
x=(\delta_*,s)\in \Omega.\] By \cite[Lemma~2.1]{cr2}, $
z_\sigma\leq u_a$ in $\Omega_{\delta_*}$, for every $\sigma\in
(0,\widetilde\delta-\delta_*)$. Letting $\sigma \searrow 0$, we
conclude Step~3.

Thus, by Steps~2 and 3, we have
\[ \xi^-\leq \liminf_{d(x)\searrow 0}\frac{u_a(x)}{h(d(x))} \leq
\limsup_{d(x)\searrow 0}\frac{u_a(x)}{h(d(x))}\leq \xi^+. \]
Taking $\ep\ri 0$, we obtain (\ref{fjk}). This finishes the proof
of Theorem~\ref{uni1}. \qed

\section{Proof of Theorem~\ref{uni3}} \label{five}
Fix $a<\lambda_{\infty,1}$ and denote by $u_a$ the unique large
solution of (\ref{0a}).

Let $\ep\in (0,1/2)$ be arbitrary and $\delta>0$ be such that (i),
(ii), (iv) from \S{\ref{nup}} are satisfied.

By \eq{tilde B} and Remark~\ref{tur}, we can diminish $\delta>0$
such that
\begin{equation}
\left\{\begin{aligned} & 1+(\widetilde c-\ep)d^\theta<b(x)/k^2(d)<
1+(\widetilde c+ \ep)d^\theta, \ \forall x\in
\Omega\ \mbox{with}\ d\in (0,\delta),\\
& k^2(t)\left[1+(\widetilde c-\ep)t^\theta\right]\ \mbox{is
increasing on}\ (0,\delta).
\end{aligned} \right.
\label{chn}
\end{equation}

Define $u^\pm(x)=\xi_0 h(d)\left[1+\chi_\ep^\pm (-\ln
d)^{-\tau}\right]$ for $x\in \Omega$ with $d\in (0,\delta)$, where
$\chi_{\ep}^\pm= \widetilde \chi \pm \ep$.

We can assume $u^\pm (x)>0$ for every $x\in\Omega$ with $d(x)\in
(0,\delta)$.

By the Lagrange mean value theorem, we obtain
\[ f(u^\pm (x))=f(\xi_0
h(d))+\xi_0 \chi_\ep^\pm \frac{h(d)} {(-\ln d)^\tau}\,f'(\Psi^\pm
(d)),\] where $ \Psi^\pm (d)=\xi_0
h(d)\left[1+\chi_\ep^\pm\lambda^\pm(d) (-\ln d)^{-\tau}\right],$
for some $\lambda^\pm (d)\in [0,1]$.

Since $f(u)/u^{\rho+1}$ is slowly varying, by Proposition~\ref{pn}
we find \begin{equation} \label{lim2} \lim_{d\searrow 0}
\frac{f(\Psi^\pm(d))}{f(\xi_0 h(d))}= \lim_{d\searrow 0}
\frac{f(u^\pm(d))}{f(\xi_0 h(d))}= 1.\end{equation}

\emph{Step} 1. There exists $\delta_1\in (0,\delta)$ so that
\begin{equation}
\left\{\begin{aligned} & \Delta u^+ +au^{+}-k^2(d)[1+(\widetilde
c-\ep)d^\theta]f(u^+) \leq 0, \quad \forall x\in \Omega \mbox{
with } d<\delta_1,\\
& \Delta u^- +au^{-}-k^2(d)[1+(\widetilde c+\ep)d^\theta]f(u^-)
\geq 0, \quad \forall x\in \Omega \mbox{ with } d<
\delta_1.\end{aligned} \right. \label{yy1}
\end{equation}

For every $x\in \Omega$ with $d\in (0,\delta)$, we have
\begin{equation} \label{min}
\Delta u^{\pm}+au^{\pm}-k^2(d)\left[1+(\widetilde c\mp\ep)
d^\theta\right]f(u^\pm)=\xi_0\, \frac{h''(d)}{(-\ln
d)^\tau}{\mathcal J}^\pm(d)
\end{equation}
where
\[\begin{aligned} {\mathcal J}^\pm(d):=& \left[
\chi_{\ep}^{\pm}\Delta d\,\frac{h'(d)}{h''(d)}+
\frac{h'(d)}{dh''(d)}\left(d(-\ln d)^\tau \Delta d-\frac{2\tau
\chi_{\ep}^{\pm}}{\ln d}\right)+
a\,\frac{h(d)}{h''(d)}\left(\chi_{\ep}^{\pm}+(-\ln
d)^\tau\right)\right.\\
&\ \, +\frac{\tau \chi_{\ep}^{\pm} \,h(d)}{d^2h''(d)\ln d}\left(1+
\frac{\tau+1}{\ln d}-d\Delta d\right)+ (-\widetilde c\pm
\ep)d^{\theta}(-\ln d)^\tau\,\frac{k^2(d)f(\xi_0 h(d))}{\xi_0
h''(d)}\\
& \ \,\left. +(-\widetilde c\pm \ep)\chi_{\ep}^{\pm}d^\theta\,
\frac{k^2(d)h(d)f'(\Psi^{\pm}(d))}{h''(d)}+{\mathcal H}(d)+
{\mathcal J}_{1}^{\pm}(d)\right].
\end{aligned} \]
Here ${\mathcal H}$ is defined by \eq{pik}, while
\[ {\mathcal J}_{1}^{\pm}(d):=
 \chi_{\ep}^{\pm}\left(1-\frac{k^2(d)h(d)f'(\Psi^{\pm}(d))}
{h''(d)}\right).\] By Lemma~\ref{aux} and \eq{lim2}, we infer that
\[\lim_{d\searrow 0}\frac{k^2(d)h(d)f'(\Psi^{\pm}(d))}{h''(d)}=
\lim_{d\searrow 0}\frac{\Psi^{\pm}(d)f'(\Psi^{\pm}(d))}
{f(\Psi^{\pm}(d))}\,\frac{k^2(d)f(\xi_0 h(d))}{\xi_0 h''(d)}
=\rho+1. \] Hence, $ \lim_{d\searrow 0} {\mathcal
J}_1^\pm(d)=-\rho \chi_\ep^\pm:=-\rho (\widetilde \chi\pm \ep)$.
Using Lemmas~\ref{aux} and \ref{goal}, we find
\[\lim_{d\searrow 0} {\mathcal J}^+(d)=-\rho\ep<0\ \
\mbox{and}\ \ \lim_{d\searrow 0}{\mathcal J}^-(d)=\rho\ep>0 .\]
Therefore, by (\ref{min}) we conclude (\ref{yy1}).

\emph{Step} 2. There exists $M^+$, $\delta^+>0$ such that
\[ u_a(x)\leq u^+(x)+M^+, \qquad \forall x\in \Omega\ \mbox{with}\
0<d<\delta^+.\] We only recover \eq{emu2}, the rest being similar
to the proof of Step 2 in Theorem~\ref{uni1}. Indeed, by
(\ref{yy1}), (\ref{chn}) and (\ref{z1}), we obtain
\[\begin{aligned}
 -\Delta u^*_\sigma (x)
&\geq a u^+(d-\sigma,s)-[1+(\widetilde c-\ep)(d-\sigma)^\theta]
k^2(d-\sigma)f(u^+(d-\sigma,s)) \\
&\geq  a u^+(d-\sigma,s)-[1+(\widetilde c-\ep)d^\theta]
k^2(d)f(u^+(d-\sigma,s))\\
& \geq  a u^+(d-\sigma,s)-b(x)f(u^+(d-\sigma,s))\\
&\geq  a (u^+(d-\sigma,s)+M^+)-b(x)f(u^+(d-\sigma,s)+M^+)\\
&= au^*_\sigma(x)-b(x)f(u^*_\sigma(x))\quad \mbox{in } {\mathcal
N}_{\sigma}.
\end{aligned}\]

\emph{Step} 3. There exists $M^-$, $\delta^->0$ such that \[
u_a(x)\geq u^-(x)-M^-, \qquad \forall x\in \Omega\ \mbox{with}\
0<d<\delta^-. \] We proceed in the same way as for proving
(\ref{nb}). To recover (\ref{em}) (with $\lambda$ given by
(\ref{la})), we show that $\Delta v_\sigma^*+a v_\sigma^*\geq
b(x)f(v_\sigma^*)$ in $\Omega_{\delta_2/4}$. Indeed, using
(\ref{chn}), (\ref{yy1}) and $(A_1)$, we find
\[\begin{aligned}
 \Delta v^*_\sigma (x)+av_\sigma^*(x)  &=
\lambda(\Delta u^-(d+\sigma,s)+au^-(d+\sigma,s))\\
 &\geq  \lambda k^2(d+\sigma)[1+(\widetilde c+\ep)(d+\sigma)^\theta]
f(u^-(d+\sigma,s)) \\
 &\geq  k^2(d)[1+(\widetilde c+\ep)d^\theta]f(\lambda u^-(d+\sigma,s))\\
 &\geq b(x)f(v^*_\sigma(x)),\quad
\forall x=(d,s)\in \Omega_{\delta_2/4}.
\end{aligned}\]

Since $\lim_{d\searrow 0}k^2(d)f(\lambda^2 u^-(x))/u^-(x)=\infty$,
there exists $\widetilde \delta\in (0,\delta_2/4)$ such that
\neweq  k^2(d)[1+(\widetilde c+\ep)d^\theta]f(\lambda^2 u^-)/u^-
\geq \lambda^2|a|,\quad \forall x\in \Omega\ \mbox{with}\ 0<d\leq
\widetilde\delta. \label{off} \endeq By Lemma~\ref{aux}, we infer
that $u^-(x)$ decreases with $d$ when $d\in (0, \widetilde
\delta)$ (if necessary, $\widetilde \delta>0$ is diminished).
Choose $\delta_*\in (0,\widetilde \delta)$ close enough to
$\widetilde\delta$ such that
\neweq
\frac{h(\delta_*)(1+\chi_\ep^-(-\ln \delta_*)^{-\tau})}
{h(\widetilde\delta)(1+\chi_\ep^-(-\ln
\widetilde\delta)^{-\tau})}<1+\lambda. \label{win}
\endeq
Hence, we regain \eq{bim}, \eq{ja} and \eq{bem}.

By (\ref{chn}) and (\ref{yy1}), we see that \eq{hy} follows if
\begin{equation} k^2(d +\sigma)[1+(\widetilde c+\ep)(d+\sigma)^\theta]
\left[f(u^-(d+\sigma,s))-f(z_\sigma(d,s))\right] \geq
a(1-\lambda)u^-(\delta_*,s)\label{ew}
\end{equation} for each $(d,s)\in
\Omega_{\delta_*}$. Using \eq{bem}, together with \eq{off} and
\eq{ja}, we arrive at \eq{ew}. From now on, the argument is the
same as before. This proves the claim of Step 3.

By Steps~2 and 3, it follows that
\begin{equation}\left\{\begin{aligned}
& \chi_\ep^+\geq \left[-1+\frac{u_a(x)}{\xi_0 h(d)}\right] (-\ln
d)^\tau -\frac{M^+(-\ln d)^\tau}{\xi_0 h(d)}\,,\ \forall x\in
\Omega \ \mbox{with}\ d<\delta^+\\
& \chi_\ep^-\leq \left[-1+\frac{u_a(x)}{\xi_0 h(d)}\right] (-\ln
d)^\tau +\frac{M^-(-\ln d)^\tau}{\xi_0 h(d)},\ \forall x\in \Omega
\ \mbox{with}\ d<\delta^-.\end{aligned}\right. \label{ulm}
\end{equation} Using Lemma~\ref{aux}, we have \[
\lim_{t\searrow 0}\frac{(-\ln t)^\tau}{h(t)}=\lim_{t\searrow 0}
\left(\frac{-\ln t}{\ln h(t)}\right)^\tau\,\frac{(\ln
h(t))^\tau}{h(t)}= \left(\frac{\rho \ell_1}{2}\right)^\tau
\,\lim_{u\to \infty} \frac{(\ln u)^\tau}{u}=0. \] Passing to the
limit $d\searrow 0$ in (\ref{ulm}), we obtain
\[\chi_\ep^- \leq \liminf_{d\searrow 0}
\left[-1+\frac{u_a(x)}{\xi_0 h(d)}\right](-\ln d)^\tau \leq
\limsup_{d\searrow 0} \left[-1+\frac{u_a(x)}{\xi_0
h(d)}\right](-\ln d)^\tau \leq \chi_\ep^+. \] By sending $\ep$ to
0, the proof of Theorem~\ref{uni3} is finished. \qed

\end{document}